\documentclass{IEEEtran4PSCC}

\ifCLASSINFOpdf
   \usepackage[pdftex]{graphicx}
\else
   \usepackage[dvips]{graphicx}
\fi

\usepackage[cmex10]{amsmath}
\usepackage{amsfonts}
\usepackage{cellspace}
\usepackage{makecell}
\setcellgapes{3.5pt}
\usepackage[font=small,labelfont=bf]{caption}
\usepackage{subcaption}
\usepackage{xcolor}
\usepackage{booktabs}
\usepackage{multirow}
\usepackage{makecell}
\usepackage{soul}
\usepackage{breqn}
\usepackage[ruled,vlined]{algorithm2e}
\usepackage{algpseudocode}

\makeatletter
\let\old@ps@headings\ps@headings
\let\old@ps@IEEEtitlepagestyle\ps@IEEEtitlepagestyle
\def\psccfooter#1{%
    \def\ps@headings{%
        \old@ps@headings%
        \def\@oddfoot{\strut\hfill#1\hfill\strut}%
        \def\@evenfoot{\strut\hfill#1\hfill\strut}%
    }%
    \def\ps@IEEEtitlepagestyle{%
        \old@ps@IEEEtitlepagestyle%
        \def\@oddfoot{\strut\hfill#1\hfill\strut}%
        \def\@evenfoot{\strut\hfill#1\hfill\strut}%
    }%
    \ps@headings%
}
\makeatother

\psccfooter{%
        \parbox{\textwidth}{\hrulefill \\ \small{23rd Power Systems Computation Conference} \hfill \begin{minipage}{0.2\textwidth}\centering \vspace*{4pt} \includegraphics[scale=0.06]{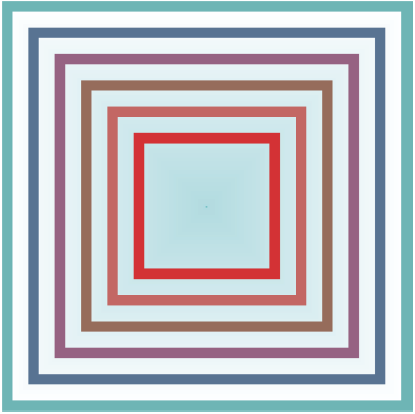}\\\small{PSCC 2024} \end{minipage} \hfill \small{Paris, France --- June 4 -- 7, 2024}}%
}

\begin{document}

\title{Predictive Optimization of Hybrid Energy Systems with Temperature Dependency}

\author{
\IEEEauthorblockN{Tanmay Mishra, Amritanshu Pandey, Mads R. Almassalkhi}
\IEEEauthorblockA{Department of Electrical and Biomedical Engineering, University of Vermont, USA\\
\{tmishra, apandey1, malmassa\}@uvm.edu}
\vspace{-0.5cm}
}

% make the title area
\maketitle

\begin{abstract}
Hybrid Energy Systems (HES), amalgamating renewable sources, energy storage, and conventional generation, have emerged as a responsive resource for providing valuable grid services. Subsequently, modeling and analysis of HES have become critical, and the quality of grid services hedges on it. Currently, most HES models are temperature-agnostic. However, temperature-dependent factors can significantly impact HES performance, necessitating advanced modeling and optimization techniques. With the inclusion of temperature-dependent models, the challenges and complexity of solving optimization problems increase. In this paper, the electro-thermal modeling of HES is discussed. Based on this model, a nonlinear predictive optimization framework is formulated. A simplified model is developed to address the challenges associated with solving nonlinear problems (NLP). Further, projection and homotopy approaches are proposed. In the homotopy method, the NLP is solved by incrementally changing the C-rating of the battery. Simulation-based analysis of the algorithms highlights the effects of different battery ratings, ambient temperatures, and energy price variations. Finally, comparative assessments with a temperature-agnostic approach illustrate the effectiveness of electro-thermal methods in optimizing HES.
\end{abstract}


\thanksto{\noindent This material is based upon work supported by the U.S. Department of Energy's Office of Energy Efficiency and Renewable Energy (EERE) under the Solar Energy Technologies Office Award Number DE-EE0010147. The views expressed herein do not necessarily represent the views of the U.S.Department of Energy or the United States Government.}
\vspace{-0.5cm}
\section{Introduction}
The imperative to decarbonize the electric power sector and fortify grid resilience has sparked heightened interest in the seamless integration of renewable energy sources \cite{GIELEN2019}. Distributed generation (DG) localized generation and enhanced grid robustness \cite{Muhtadi2021}. However, the inherent intermittency of renewable sources poses significant challenges to grid stability and reliability. Addressing these limitations requires a transformative shift towards adopting cutting-edge Hybrid Energy Systems (HES), intelligently combining diverse energy sources and advanced storage technologies \cite{Roy2022}. Through precise orchestration of energy source integration, these HES exhibit remarkable potential in mitigating the drawbacks associated with DG and microgrids, thereby ensuring reliable power \cite{Murphy2021}. Through a cohesive approach that harnesses renewable energy sources alongside state-of-the-art energy storage, HES offers enhanced grid stability, improved energy management, increased resilience, flexibility, reduced emissions and cost-effectiveness through optimized energy management, leveraging diverse energy sources and storage technologies with greater energy security over DG and micro-grids \cite{Zhang2019,LI2018}.

As the number of HES technologies increases, the technology gap associated with their analysis becomes more significant \cite{Badwawi2015}. This significance is accentuated under extreme weather conditions, where the temperature dependence of various HES components becomes crucial. The 2021 winter storm Uri, Texas proved a fantastic revenue opportunity for batteries, enabling some systems to return multiples of their capital installation costs in a single year. The battery resources with the highest and lowest revenue differed by two primary factors: bidding optimization and time in the market \cite{Gridmatic2022}. Similarly, the September 2022 California heat wave highlighted batteries' substantial role by contributing a significant portion of their capacity to the market, impacting regulation dynamics. Operating constraints, particularly those tied to the state of charge (SoC) and related factors, emerge as pivotal under exigent circumstances like extreme weather \cite{CAISO2023}. These variables collectively underscore the disproportionate revenue impact experienced amidst extreme weather conditions.

Many battery manufacturers ignore ambient weather's influence on performance, opting for large, insulated enclosures to mitigate environmental effects using HVAC systems. As demand for fast-charging batteries rises, HVAC complexity grows. Consequently, battery dispatch scheduling becomes highly temperature-dependent and must consider HVAC performance and power consumption. To enable predictive scheduling in such systems, temperature-dependent models are crucial. These models, tailored to high-rate discharge batteries, ensure efficient utilization amidst varying weather conditions, aligning battery performance with evolving energy system needs~\cite{TOMASZEWSKA2019}. Optimal scheduling of HES using temperature-dependent models entails several challenges, such as model complexity, the introduction of non-linearities into the system's behavior, mixed integer nature due to discrete operating modes, increased dynamic nature, and limited data associated with thermal variables.

Numerous contemporary studies have explored solutions to these problems. In \cite{GUANGQIAN2018}, a machine learning-based algorithm was proposed to optimize the HES. However, this approach does not prevent the possibility of simultaneous charging and discharging (SCD) in the battery.
In~\cite{garifi2019}, the emergence of suboptimal solutions due to SCD has been discussed and a convex battery model is introduced to ensure optimal outcomes. 
Notably, this model is applicable only if battery power is not exported back to the grid. In~\cite{Nazir2020}, a weighted battery power term was introduced into the objective function to mitigate SCD in multi-period distributed energy resources (DER) dispatch, aiming at objectives such as voltage regulation, line loss minimization, and power reference point tracking. A non-linear IVQ model for battery optimization is proposed in \cite{Aaslid2020}, where a cubic spline (a piecewise polynomial function) is fitted to experimental data from a Nissan Leaf cell battery. This process aims to model losses and establish relationships between voltage, current, and charge. However, since this model is specifically tailored to the characteristics of a particular battery, its applicability to different battery types remains uncertain. Most of these models do not account for thermal constraints limiting battery charge and discharge. 
Temperature effect and thermal impact on the performance of Li-ion batteries are detailed in \cite{Ma2018}. Hence, a comprehensive electro-thermal model was developed that emphasizes the physical attributes of the battery and associated losses \cite{Schimpe2018}. Different electro-thermal battery models are presented and it can be observed that battery system complexity increases with thermodynamic considerations~\cite{Rosewater2019}. 

This paper discusses a temperature-dependent model for HES based on physical dynamics. It formulates and analyzes a non-convex, non-linear problem (NLP) for predictive HES scheduling using temperature-dependent models. It introduces a simplified mixed-integer problem to obtain the best solution under simplified constraints while identifying the bounds of the objective. However, this simplified solution may not guarantee feasibility in the original NLP. Therefore, a projection algorithm is proposed to find the closest feasible NLP solution to the best solution obtained from the simplified model. Additionally, an incremental homotopy algorithm is introduced to ensure the best NLP solution over the battery's C-ratings. The paper presents various case studies demonstrating the performance of both the NLP and simplified HES models, along with the projection and homotopy algorithms. Furthermore, a temperature-agnostic HES model has been used and compared with the original NLP model's solution to assess its dispatch schedule's feasibility. The predictive optimization has been stress-tested against the negative energy prices. The PV curtailment approach has been integrated into problem formulation to avoid the loss of revenue during negative pricing. 

The structure of the paper is given as: Section~\ref{system_des} and Section~\ref{Model}, detail HES and its subcomponent models including temperature dependencies, respectively. Section~\ref{Problem formulation} formulates a revenue maximization problem and discusses associated challenges. Section~\ref{case_studies} presents performance comparisons of HES models under different conditions. Section~\ref{Conclusion} provides the paper's conclusion.

\section{Hybrid Energy System Description}\label{system_des}

A HES represents an amalgamation of energy generation, storage, and conversion assets whose coordination achieves desired operational objectives. For example, HES applications can range from GW-scale multi-energy systems that leverage waste heat from thermal generators to improve electrolyzer production to electricity-only, colocated battery+PV off-grid microgrids to multi-locational virtual power plants (VPPs)~\cite{Murphy2021}. These hybrid PV+battery systems have the potential to reduce costs and increase the energy output compared to separate PV and battery storage systems of similar size. Today, PV+battery HES are being deployed at an increasing rate and is the focus of this paper~\cite{Eurek2021}.

\begin{figure}[!htb]
\centering
\includegraphics[width=0.38\textwidth]{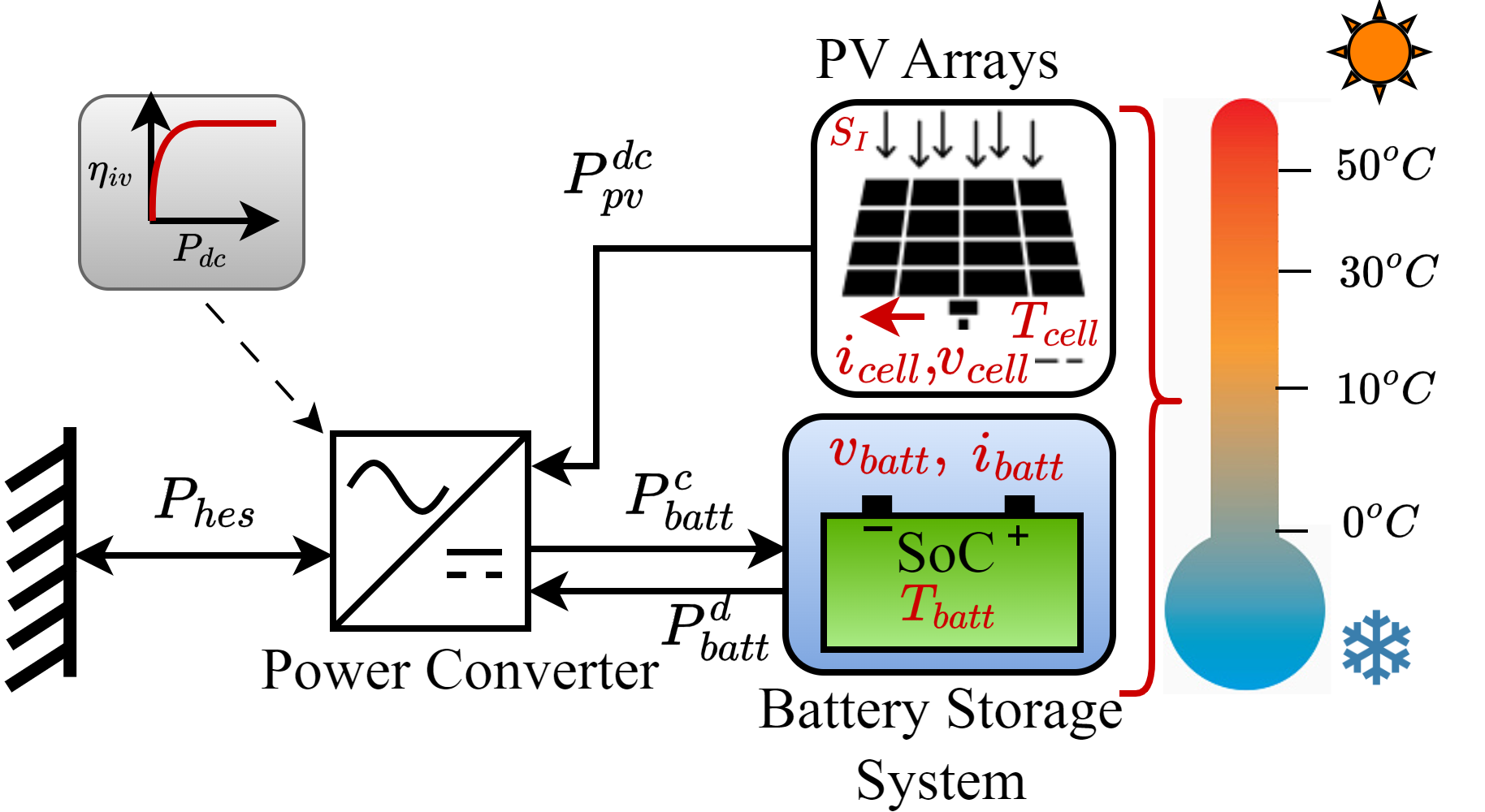}
\caption{Schematic diagram of hybrid energy system.}
\label{HES_schematic}
\centering
\end{figure}

The schematic of a PV+ battery storage HES is shown in Fig.~\ref{HES_schematic}. The PV and battery storage systems are tied to the grid through a power converter having part-load efficiency. The part-load efficiency of the power converter depends on the input DC power. The total power injected into the grid, denoted as $P_{\text{hes}}$, is calculated as the combined sum of PV and battery power contributions, with the ambient temperature dependency. In the figure, $P^{\text{c}}_{\text{batt}}$ and $P^{\text{d}}_{\text{batt}}$ represent charging and discharging power respectively, while $P^{\text{dc}}_{\text{pv}}$ is the total dc power output of the PV arrays with $i_{\text{cell}}$ and $v_{\text{cell}}$ as dc currents and voltages of each PV cell, respectively. Further details regarding the electro-thermal battery model can be found in the subsequent subsection as shown in Fig.~\ref{Battery_Thermal_ECM}. 
\section{Temperature-dependent HES Models}\label{Model}
The paper focuses on optimizing accurate temperature-dependent models for HES constituents, including PV, battery storage, and inverters. These models pave the way for meticulous fine-tuning of energy generation, storage, and distribution, thereby orchestrating effective energy management and informed decision-making. The ensuing sections delve into the temperature-dependent modeling intricacies of the HES sub-components.

\subsection{Power converter model}
The power converter interfaces the HES sub-components (PV and battery) with the grid, assuming a stiff grid voltage. The relationship between the grid-side power $P^{\text{ac}}$ and the input dc power $P^{\text{dc}}$ is determined by the part-load inverter efficiency $\eta_{\text{iv}}$ and it is given by~\cite{Gilman2015}:
\begin{subequations}
\begin{align}
P^{\text{ac}}(t)&= \eta_{\text{iv}}(t) P^{\text{dc}}(t)\\
\eta_{\text{iv}}(t)&= \frac{\eta_{\text{max}}}{1+e^{-\gamma P^{\text{dc}}(t)}},
\end{align}
\label{Inv_eff}
\end{subequations}
where $\gamma$ represents the gradient of the inverter's efficiency curve. The DC side of the power converter is connected to the PV and battery. It has been observed from the test results presented in \cite{Sandia2007} that the inverter efficiencies typically do not have a strong temperature dependence. Subsequently, the PV and battery systems connected to the DC side of the power converter are discussed.

\subsection{The Photovoltaic (PV) model}
One of the simplest and most popular PV model for system-level studies, is based on the efficiency of the PV cells ($\eta_{\text{cell}}$) and solar irradiance ($S_{\text{I}}(t)$), which is given by~\cite{Gilman2015}:
\begin{equation}
 P_{\text{cell}}(t) = \eta_{\text{cell}}A_{\text{pv}}S_{\text{I}}(t)  
 \label{P_cell_temp_agn}
\end{equation}
where $P_{\text{cell}}(t)$ represents the power generation from each cell, and $A_{\text{pv}}$ is the area of PV cells in $m^2$. However, this model has limitations, notably due to the exclusion of semiconductor physics, temperature impacts, and contact losses. It has been used as a \textbf{temperature-agnostic PV model} in this paper. This underscores the necessity for a more accurate PV characterization model. These limitations can be addressed by building the PV model from the cell level and expanding it to PV arrays. For this work, an equivalent circuit model of PV has been used \cite{Masters2004}, whose behavior can be described by:
\begin{equation}
    i_{\text{cell}}(t) = i_{\text{sc}}(t)-I_{\text{0}} \left(e^{38.9v^{\text{pv}}_{\text{oc}}(t)}-1\right)-\left(\frac{v^{\text{pv}}_{\text{oc}}(t)}{R_{\text{P}}}\right)
    \label{PV_current}
\end{equation}
\noindent where $i_{\text{sc}}(t)$ is the short circuit current and $v^{\text{pv}}_{\text{oc}}(t)$ is the open circuit terminal voltage, i.e., the photodiode voltage when conducting. The constant value~38.9 is calculated at 300K operating temperature from the charge and Boltzmann's constant. $R_{\text{P}}$ is the parallel resistance denoting the shading effect. Here, the $i_{\text{sc}}(t)$ is proportional to the irradiance $S_{\text{I}}(t)$. The voltage and power equations are given by:
\begin{subequations}
\begin{align}
v_{\text{cell}}(t) &= v^{\text{pv}}_{\text{oc}}(t)-i_{\text{cell}}(t)R_{\text{S}}\\
P_{\text{cell}}(t) &= i_{\text{cell}}(t)v_{\text{cell}}(t)\\
P^{\text{g}}_{\text{pv}}(t) &=  ({L~M~K})~P_{\text{cell}}(t)
\end{align}
\label{PV_Cell_Power}
\end{subequations}
\noindent where $R_{\text{S}}$ is the series resistance denoting the contact losses, $L$ is the series connected modules per string, $M$ is parallel strings and $K$ is the number of PV arrays. The PV output power is sensitive to changes in temperature. The relation between ambient temperature  ($T_{\text{amb}}$) and PV cell temperature ($T_{\text{cell}}$) is given by:
\begin{equation}
    T_{\text{cell}}(t) = T_{\text{amb}}(t)+ \left(\frac{(\text{NOCT}-20^{o})}{800}\right)S_{\text{I}}(t)
    \label{Ambient_Temperature}
\end{equation}
\noindent where NOCT is nominal operating cell temperature when $T_{\text{amb}} = 20^{\circ}$C and wind speed is 1 m/s. The temperature-dependent power and voltage equations of PV cell are given by~\cite{Masters2004}:
\begin{subequations}
\begin{alignat}{3}
    v^{\text{pv}}_{\text{oc}}(t) = V^{\text{std}}_{\text{oc}}[1-C_{\text{VT}}(T_{\text{cell}}(t)-25)]\\
    P_{\text{cell}}(t) = i_{\text{cell}}(t)v_{\text{cell}}(t)[1-C_{\text{PT}}(T_{\text{cell}}(t)-25)]
\end{alignat}
\label{PV_Power_Temperature}
\end{subequations}
\noindent Here, $C_{\text{VT}}$ and $C_{\text{PT}}$ are temperature sensitivity coefficients for voltage and power respectively. The $V^{\text{std}}_{\text{oc}}$ is the open circuit voltage at standard temperature of $T_{\text{cell}} = 25^{o}C$. 
The described PV model is nonlinear due to the diode current and the bilinear $i-v$ relation to cell power. However, solar PV output is also inherently variable, which can strain grid stability. When faced with high PV generation, low demand, or negative energy prices, PV curtailment can be employed as a solution~\cite{Eric2019}. PV curtailment can be achieved using a curtailment factor $\beta_{\text{cur}} \in [0,1]$. For the total PV power generated $P^{\text{g}}_{\text{pv}}$, the power exported to the grid would be denoted as $P_{\text{pv}}$. Mathematically, this relationship is expressed as follows:
\begin{equation}
    P_{\text{pv}}(t) = \beta_{\text{cur}}P^{\text{g}}_{\text{pv}}(t)
    \label{Curtailment}
\end{equation}
Secondly, By adding batteries to PV installations to engender a HES, energy can be stored during periods of excess generation and discharged when needed. 
Thus, (partly) overcoming solar PV's intermittency, thereby enhancing the reliability, utilization, and sustainability of HES.

\subsection{The battery model}
\begin{figure}[!htb]
\centering
\includegraphics[width=0.38\textwidth]{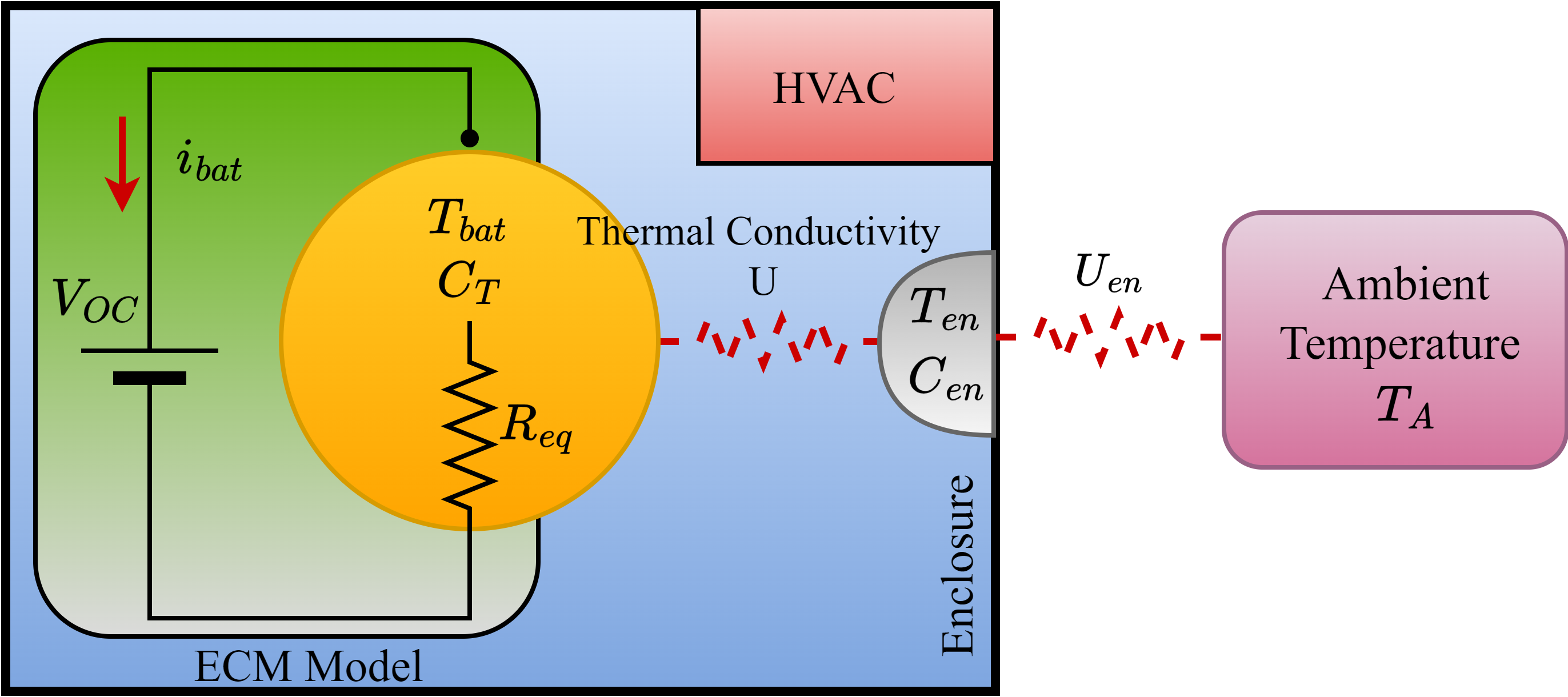}
\caption{Battery's ECM with enclosure and HVAC.}
\label{Battery_Thermal_ECM}
\end{figure}

In the realm of battery modeling, a spectrum of approaches graces the literature, encompassing both rudimentary approximations and intricate physical precision. In the context of HES applications, an ideal pursuit is a model that seamlessly merges computational efficiency with versatility~\cite{Sioshansi2019}. The Energy Reservoir Model (ERM) is commonly used in system level research. In this model, the battery's SoC dynamics, $E(t) \in [\underline{E}, \bar{E}]$, depend on the charging and discharging power, $P^{\text{c}}_{\text{batt}}$ and $P^{\text{d}}_{\text{batt}}$, respectively. The ERM model is as follows~\cite{Ferreira2019}:
\begin{subequations}
\begin{equation}
E(t)=E_{\text{0}} + \frac{1}{E_{\text{c}}}\int_{0}^{t} \left(\eta^{\text{batt}}_{\text{iv}}P^{\text{c}}_{\text{batt}}(\tau)
+ \frac{1}{\eta^{\text{batt}}_{\text{iv}}}P^{\text{d}}_{\text{batt}}(\tau)\right) d\tau
\end{equation}
\begin{equation}
\text{SoC limits}:~ \underline{E} \le E(t) \le \bar{E} 
\end{equation}
\begin{equation}
0 \le P^{\text{c}}_{\text{batt}}(t) \le P^{\max}_{\text{batt}}, ~-P^{\max}_{\text{batt}}\le P^{\text{d}}_{\text{batt}}(t) \le 0
\end{equation}
\label{ERM_model}  
\end{subequations}
\noindent where $E_{\text{c}}$ is the energy capacity of the battery in $kWh$. This model is sufficiently accurate when the battery operates over a small voltage range but has inaccuracies when a larger range of voltage values are encountered~\cite{Rosewater2019}. Further, this model ignores temperature's impact on SoC and charging/discharging schedules. Striking a balance between exactitude and intricacy, an electro-thermal Equivalent Circuit Model (ECM) for Li-ion battery is used as shown in Fig.~\ref{Battery_Thermal_ECM}. A zeroth-order ECM (second-order model with neglected dynamics) has been used. As our study focuses on steady-state systems, we have considered only resistive elements~\cite{Morstyn2020}.

Mathematically, this model can be defined by:
\begin{subequations}
\begin{align}
C_{\text{Ah}}\frac{\partial E(t)}{\partial t} &= i_{\text{batt}}(t)\\
v_{\text{batt}}(t)&=V_{\text{OC}}+i_{\text{batt}}(t)R_{\text{eq}}\\
V_{\text{OC}} &= V_{\text{m}} E(t)+ V_{\text{0}} 
\end{align}
\label{battery_model}
\end{subequations}
\noindent where $i_{\text{batt}}$ and $v_{\text{batt}}$ are battery net current and terminal volatge respectively, $R_{\text{eq}}$ is equivalent internal resistance, $C_{\text{Ah}}$ is charge capacity in \textit{Ah} and $V_{\text{OC}}$ is open circuit battery voltage at $100\%$ SoC. The slope and intercepts of a linear $V_{\text{OC}}$ model are represented by $V_{\text{m}}$ and $V_{\text{0}}$, respectively. Maintaining battery temperature within limits is crucial for preserving SoC and state of health (SoH), preventing damage and safety risks, particularly in hot conditions or high-power usage.
% The battery temperature can significantly affect SoC and state of health (SoH), especially in hot environments or in high-power operations, where excessive losses drive up cell temperatures. Thus, enforcing temperature constraints is essential to prevent adverse effects like thermal run-away, which can damage the battery and pose safety risks. 
To extend the thermal capability of a battery, an HVAC system, as shown in Fig. \ref{Battery_Thermal_ECM}, can be used to keep the battery temperature within its design limits i.e. $T_{\text{batt}} \in [T^{\min}_{\text{batt}}, T^{\text{max}}_{\text{batt}}]$. Considering battery within an enclosure with an HVAC system, overall thermal dynamics are given by the following two equations:
\begin{subequations}
\begin{alignat}{3}
    C_{T}\frac{\partial T_{\text{batt}}(t)}{\partial t} &= R_{\text{eq}}i^{2}_{\text{batt}}(t)+ U(T_{\text{en}}(t)-T_{\text{batt}}(t)) \\
    C_{\text{en}}\frac{\partial T_{\text{en}}(t)}{\partial t} &= N_{\text{cell}}\left[U(T_{\text{batt}}(t)-T_{\text{en}}(t))\right]+U_{\text{en}}(T_{\text{amb}}(t)\\
    & -T_{\text{en}}(t))-\eta_{\text{hv}}P_{\text{hv}}(t) \notag
\end{alignat}
\end{subequations}
where $C_{\text{T}}$ and $C_{\text{en}}$  are the heat capacity of the lumped volume and battery enclosure respectively, \textit{U} represents the thermal transmittance between the battery surface and the enclosure, while $U_{\text{en}}$ pertains to the thermal transmittance between the enclosure and the surrounding environment, and $N_{\text{cell}}$ is the number of battery cells in the enclosure. The $\eta_{\text{hv}}$ and $P_{\text{hv}}$ are efficiency and power drawn by HVAC, respectively. This model assumes constant or no airflow. The battery power is given as follows:
\begin{subequations}
\begin{align}
P^{\text{c}}_{\text{batt}}(t) & = i^{\text{c}}_{\text{batt}}(t)v_{\text{batt}}(t)\\
P^{\text{d}}_{\text{batt}}(t) & = i^{\text{d}}_{\text{batt}}(t)v_{\text{batt}}(t)\\
i_{\text{batt}}(t) & = i^{\text{c}}_{\text{batt}}(t)+i^{\text{d}}_{\text{batt}}(t),
% P_{\text{hv}}(t)       & = P^{\text{cool}}_{\text{hv}}(t)+P^{\text{heat}}_{\text{hv}}(t),Here
\end{align}
\label{battery_Power}
\end{subequations}
\noindent where $i^{\text{c}}_{\text{batt}}$ and $i^{\text{d}}_{\text{batt}}$ represent charging and discharging battery currents respectively with $i^{\text{c}}_{\text{batt}}(t)\ge0$ and $i^{\text{d}}_{\text{batt}}(t)\le0$. Using these models of the inverter, PV and battery, the predictive optimization problem has been formulated and discussed in the next section.

\section{Predictive Optimization of HES}\label{Problem formulation}
The predictive optimization enables the seamless integration of various energy sources. It extends to providing grid services (such as voltage regulation and frequency control) and efficient use of energy storage by predicting energy demand and supply fluctuations. This includes optimizing charge and discharge cycles to reduce costs and elongate the system lifespan. In this paper, predictive optimal scheduling of the HES is formulated  as an optimization problem, subject to equality and inequality constraints governed by physics-based models of PV, battery and inverter. The objective of the optimization problem is maximization of revenue generated by the grid-connected HES. This is defined by the sum of the product of net HES power, energy price and time step. The energy price in $\$/MWh$ denoted by $P^{r}[k]$, with~$ k\in\{1\hdots N\}$. The equality constraints are defined using \eqref{Inv_eff} and \eqref{battery_model}-\eqref{battery_Power}. The inequality constraints on SoC, charging-discharging currents, battery temperature and voltage are defined as the bounds based on battery operation and performance. The decision variables are denoted by \textbf{x} =$\{i^{c}_{\text{batt}},i^{d}_{\text{batt}},P_{\text{hv}}\}$, each decision variable can be denoted by vector \textbf{$x_{i}$} where $i\in\{1\hdots 3\}$ and $\textbf{$x_{i}$} \in \mathbb{R}^{N}$.

The HES predictive optimization, rooted in physics-based models, addresses a single-objective, multi-period deterministic problem. Within this framework, two temperature-dependent models, namely NLP and simplified (MIP), alongside a temperature-agnostic model based on ERM battery model and PV model from \eqref{P_cell_temp_agn} have been formulated. Furthermore, a projection algorithm has been introduced to yield the optimal solution from the MIP problem within the NLP feasible set. Additionally, a homotopy algorithm has been proposed to guarantee the attainment of the best solution for the original NLP.
\vspace{-0.1cm}
\subsection{NLP problem}
Based on the HES model discussed in section \ref{Model}, a non-convex and non-linear (NLP) predictive optimization problem has been formulated as follows:\\
%\fbox{\begin{minipage}{23.5em}
\textbf{Objective:}
\begin{equation}
\max_{\mathbf{x}}\qquad T_{s} \sum_{k=1}^{N}P^{r}[k]P_{hes}[k]
\end{equation}
\textbf{Subject to:}
\begin{subequations}
\begin{align}
P_{\text{hes}}[k] - P_{\text{pv}}[k] + P_{\text{batt}}[k] + P_{hv}[k] = 0 \label{hes_power}\\
P_{\text{batt}}[k] - (1/\eta^{\text{batt}}_{\text{iv}}[k])P^{\text{c}}_{\text{batt}}[k] -\eta^{\text{batt}}_{\text{iv}}[k]P^{\text{d}}_{\text{batt}}[k] = 0\\
\eta^{\text{batt}}_{\text{iv}}[k]- \frac{\eta_{\text{max}}}{1+e^{-\gamma P_{\text{batt}}[k]}} = 0\label{batt_eff}\\
P^{\text{c}}_{\text{batt}}[k] - w_{\text{kw}}~(i^{c}_{\text{batt}}[k] v_{\text{batt}}[k]) = 0\label{PC}\\
P^{\text{d}}_{\text{batt}}[k] - w_{\text{kw}}~(i^{d}_{\text{batt}}[k] v_{\text{batt}}[k]) = 0 \label{PD}\\
i_{\text{batt}}[k] - i^{c}_{\text{batt}}[k]  -i^{d}_{\text{batt}}[k] =0\\
E[k] - E[k-1] - \frac{1}{C_{\text{Ah}}}i_{\text{batt}}[k-1] T_{\text{s}}=0\\
i^{\text{c}}_{\text{batt}}[k]~ i^{\text{d}}_{\text{batt}}[k] =0 \label{CS}\\
v_{\text{batt}}[k] -V_{\text{m}} E[k]-V_{\text{0}} - i_{\text{batt}}[k]R_{\text{eq}}=0\label{VSOC}\\
T_{\text{batt}}[k] - T_{\text{batt}}[k-1] - \frac{1}{C_{\text{T}}}(R_{\text{eq}}i^{2}_{\text{batt}}[k-1] \qquad \notag \\
 + U(T_{\text{en}}[k-1] - T_{\text{batt}}[k-1])) T_{\text{s}} = 0 \label{Temp1}\\
T_{\text{en}}[k] - T_{\text{en}}[k-1] - \frac{1}{C_{\text{en}}}(N_{\text{cell}}U (T_{\text{batt}}[k-1]\qquad \notag\\
 -T_{\text{en}}[k-1]) + U_{\text{en}}(T_{\text{amb}}[k-1]-T_{\text{en}}[k-1]) \qquad  \notag\\
 -\eta_{\text{hv}}P_{\text{hv}}[k-1]) T_s = 0\label{Temp2}\\
\underline{E} \le E[k] \le \bar {E} \label{SOC_ineq}\\
\underline {v_{\text{batt}}} \le v_{\text{batt}}[k] \le \bar{v_{\text{batt}}}\\
T^{\text{min}}_{\text{batt}} \le T_{\text{batt}}[k] \le T^{\text{max}}_{\text{batt}}\\
T^{\text{min}}_{\text{en}} \le T_{\text{en}}[k] \le T^{\text{max}}_{\text{en}}\\
0 \le P_{\text{hv}}[k] \le \bar{P_{\text{hv}}}\\
0 \le i^{\text{c}}_{\text{batt}}[k] \le i^{\text{max}}_{\text{batt}} \\
-i^{\text{max}}_{\text{batt}} \le i^{\text{d}}_{\text{batt}}[k] \le 0 \label{batt_current_ineq}\\
T_{\text{en}}[0] = T^{0}_{\text{en}}, T_{\text{batt}}[0] = T^{0}_{\text{batt}}\\ 
E[0]  = E_{0}, E[N] = E_{0} \label{sustain} 
\end{align}
\label{NLP}
\end{subequations}

\noindent where $w_{\text{kw}}$~=~0.001 is Watt to kW conversion factor. The complementary slackness given by \eqref{CS} avoids the simultaneous charging and discharging of the battery. The \eqref{SOC_ineq}-\eqref{batt_current_ineq} represent the inequality associated with the NLP. The \eqref{sustain} ensures sustainability, i.e. the SoC should return to the initial level at the end of the day. 

In this HES system, the PV subsystem plays no role in decision-making (except for curtailment during periods of negative energy prices); instead, it depends on inputs $S_{\text{I}}$ and $T_{\text{amb}}$. The power generated by the PV subsystem is calculated using a temperature-dependent model described by \eqref{PV_NLP}. In this case, the contact losses and shading effect in PV have been neglected. This calculated PV power is then directly used to compute the total HES power ($P_{\text{hes}}$) in \eqref{hes_power}. The PV model is given as:
\begin{subequations} \label{PV_model}
\begin{align}
P_{\text{pv}}[k]- \eta^{\text{pv}}_{\text{iv}}[k]P^{\text{dc}}_{\text{pv}}[k]=0 \label{PV1}\\
\eta^{\text{pv}}_{\text{iv}}[k]- \frac{\eta_{\text{max}}}{1+e^{-\gamma P^{\text{dc}}_{\text{pv}}[k]]}} = 0\\
P^{\text{dc}}_{\text{pv}}[k]-(\text{L~M~K}) P_{\text{cell}}[k]=0\\
P_{\text{cell}}[k]-i_{\text{cell}}[k]v_{\text{cell}}[k][1-C_{\text{PT}}(T_{\text{cell}}[k]-25)]=0\\
i_{\text{cell}}[k] -\left(\frac{i_{\text{sc}}}{S^{\text{rated}}_{\text{I}}}\right)S_{\text{I}}[k]+I_{\text{0}} \left(e^{38.9V_{\text{oc}}^{\text{pv}}[k]}-1\right)=0\\
v^{\text{pv}}_{\text{oc}}[k] -V^{\text{std}}_{\text{oc}}[1-C_{\text{VT}}(T_{\text{cell}}[k]-25)]=0\\
T_{\text{cell}}[k] - T_{\text{amb}}[k]- S_{\text{I}}[k]\left(\frac{\text{NOCT}-20^{o}}{800}\right)=0. \label{PV7}
\end{align}
\label{PV_NLP}
\end{subequations}
\noindent The above formulated NLP pose challenges in optimization due to the complexity of finding the best solution amid multiple potential optima. Secondly, the NLP problem can be very sensitive to the initial conditions.

\subsection{Simplified Model}
A significant challenge in solving the NLP problem is the prevalence of numerous local optima, often causing optimization algorithms to converge to suboptimal solutions. In the previously formulated optimization problem, the equations \eqref{batt_eff}-\eqref{PD}, \eqref{CS}-\eqref{Temp1} introduce non-linearity and non-convexity. Specifically, equations \eqref{PC} and \eqref{PD} exhibit bilinear characteristics. To mitigate this bilinearity effect, it is possible to make an approximation by considering a
constant battery voltage $\overline{V}_{\text{batt}}$ at a specific SoC. The simplified equations are given by: 
\begin{subequations}
\begin{align}
P^{\text{c}}_{\text{batt}}[k] - i^{\text{c}}_{\text{batt}}[k] \overline{V}_{\text{batt}} = 0\\
P^{\text{d}}_{\text{batt}}[k] - i^{\text{d}}_{\text{batt}}[k] \overline{V}_{\text{batt}} = 0 
\end{align}
\label{PC_PD_simplified}
\end{subequations}
\noindent for $k\in\{1\hdots N\}$. This approach linearizes the bilinearity associated with the $P^{\text{c}}_{\text{batt}}$ and $P^{\text{d}}_{\text{batt}}$ w.r.t the $v_{\text{batt}}$. Secondly, The non-convexity due to the complementary slackness condition can be reformulated as a mixed-integer problem (MIP), presented as follows \cite{garifi2019}:
\begin{subequations}
\begin{align}
0\le i^{\text{c}}_{\text{batt}}[k] \le Z[k]i^{\text{max}}_{\text{batt}}\\ 
(Z[k]-1)i^{\text{max}}_{\text{batt}} \le i^{\text{d}}_{\text{batt}}[k] \le 0
\end{align}
\label{MIP}
\end{subequations}
\noindent for $k\in\{1\hdots N\}$. Where $ Z[k]\in \mathbb \{0,1\}$ such that if Z[k] = 1 $\rightarrow$ $i^{\text{d}}_{\text{batt}}[k]$ = 0 and $i^{\text{c}}_{\text{batt}}[k]\in \mathbb [0,i^{\text{max}}_{\text{batt}}]$. 
To handle non-linearity in \eqref{Temp1}, we apply second-order conic programming (SOCP) relaxation, which linearizes quadratic terms for efficient optimization while preserving convexity. The SOCP relaxation is expressed as follows:
\begin{align}
\alpha[k] \ge i^{2}_{\text{batt}}[k]
\label{alpha}
\end{align}
for $k\in\{1\hdots N\}$. The non-linearity associated with the inverter efficiency is addressed by taking $\eta_{\text{batt}}$ = 0.95.  Now, the original NLP has been reformulated as mixed integer linear program using \eqref{PC_PD_simplified}-\eqref{alpha}. The robust MIP solvers (e.g. Gurobi, CPLEX, etc) enable feasible and fast solutions leveraging the optimality gap as a measure of the quality of the solution. Utilizing the solution derived from the relaxed problem as an initial point significantly enhances the efficiency of the NLP solver when seeking the optimal solution to the original, more intricate NLP.

\subsection{Projection}
The solutions obtained from the simplified model are most likely a global optimal based on the optimality gap. However, it may yield solutions that partially adhere to the NLP's constraints, which are based on the physics of the HES model. Hence, the projection technique is used to ensure that the obtained solutions are best and feasible, bridging the gap between the mathematical rigour and the physics of the model. The objective of the projection problem is often defined as finding the solutions that minimize the 2-norm between the solution from simplified model ($x^{*}$) and projected solutions ($x$). It is given as:\\
\textbf{Objective:}
\begin{equation}
\min_{\mathbf{x}}\| (x[k]-x^*[k])\|_2
\end{equation}
for $k\in\{1\hdots N\}$. The equality and inequality constraints are same as the NLP problem defined by \eqref{NLP}. Even though projection is still an NLP problem, it gives the solution in the vicinity of the best solution obtained from the simplified model while satisfying the feasibility of the NLP model.

\subsection{Homotopy}
\begin{figure}[!htb]
\centering
\includegraphics[width=0.4\textwidth]{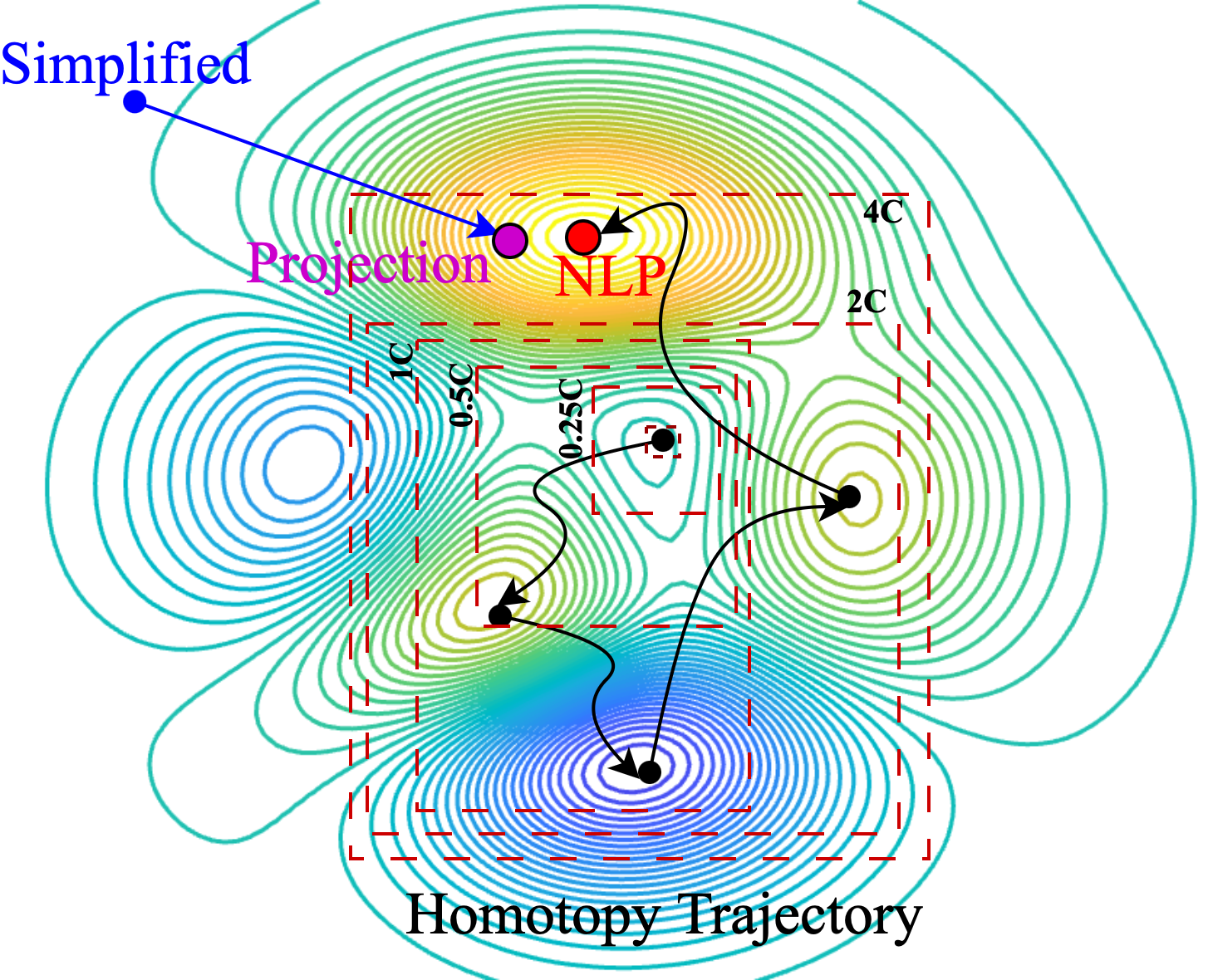}
\caption{Illustrating the different methods presented herein, including the iterative homotopy approach based on the battery's C-rating.}
\label{Contour}
\end{figure}

The homotopy algorithm is frequently used to tackle intricate non-convex optimization problems across various fields \cite{Park2020}. Homotopy employs a continuous trajectory in NLP, transitioning from the initial complex problem to a tractable one. By varying the homotopy factor (in this case, the C-rating of the battery), the problem is transformed into a sequence of sub-problems, each solved sequentially. The first subproblem is trivial, whose solution can be easily obtained, and subsequent ones converge rapidly, starting from the prior subproblem's solution, eventually leading to the original problem and facilitating robust mathematical convergence \cite{Pandey2021}. In this subsection, the homotopy method is explored to tackle battery optimization problems with varying C-ratings which is an essential factor affecting battery performance, determining the $i^{\text{max}}_{\text{batt}}$.
\\
In this case, the NLP trajectory starts with a battery C-rating of $\approx 0.0$, where the solution space is notably limited but easily navigable. From there, it progressively increases the battery's C-rating, moving towards higher C-ratings, with a small step size denoted as $\Delta C_{\text{rating}}$. This gradual transition allows the optimization algorithm to explore solutions incrementally, avoiding abrupt jumps which might lead to convergence issues. Moving through intermediate C-ratings efficiently traverses the entire spectrum of feasible solutions, ensuring potential optimum haven't been overlooked. The gradual transition from a small C-rating to higher values enables thorough exploration of the solution space, significantly enhancing the chances of finding a global solution. Fig.~\ref{Contour} illustrates the contour path depicting the NLP's feasible space and the progressive navigation through the homotopy solution space, with an emphasis on ascending C-ratings. Additionally, it provides a clear depiction of the projection of the simplified solution within the NLP's feasible space.

\vspace{0.5cm}
\begin{algorithm}
\caption{Homotopy via battery's C-rating}
\begin{algorithmic}[1]
\State \textbf{Parse:} Battery data, energy price, $T_{\text{amb}}$, $P_{\text{pv}}$ and $T_{\text{s}}$.
\State \textbf{Decision:}  $i^{\text{c}}_{\text{batt}}[k]$,  $i^{\text{d}}_{\text{batt}}[k]$, and $P_{\text{hv}}[k]$
\State \textbf{Warm start:} Enable $\verb|warm_start_init_point|$
\State \textbf{Initialize:} set $C_{\text{rating}}\approx 0$,
calculate current limits $i^{\text{max}}_{\text{batt}} := C_{\text{rating}}\times C_{\text{Ah}}$
\State \textbf{while} $C_{\text{rating}} \le  C_{\text{max}}$ %$C_{0} \rightarrow C_{\text{max}}$}
    \indent
        \begin{itemize}
            \item Solve the NLP$(C_{\text{rating}}) \Rightarrow (x^*, \lambda^*, \mu^*)$, 
            \item If primal and dual feasibility
            \item Then increment $C_{\text{rating}} \leftarrow C_{\text{rating}}  + \Delta C_{\text{rating}}$
            \item Update the previous solution to warm start NLP.
        \end{itemize}
    %\endindent
\State \textbf{end} for loop
\end{algorithmic}
\end{algorithm}

\subsection{Temperature agnostic model}\label{temp_agnotic_section}
A comparison with the temperature-agnostic model is crucial to assess the impacts and benefits of temperature-dependent model. In this subsection, the predictive optimization employs the temperature-agnostic model which is formulated using \eqref{P_cell_temp_agn} and \eqref{ERM_model} as discussed in Section~\ref{Model}. The problem formulation is as follows:\\
\textbf{Objective:}
\begin{equation}
\max_{\mathbf{x_}}\sum_{k=1}^{N}P^{r}[k]P_{hes}[k]T_{s}
\end{equation}
\textbf{Subject to:}
%Here, $N=96$ with $T_{s}=0.25$ hour. 
\begin{subequations}
\begin{align}
P_{\text{hes}}[k] - P_{\text{pv}}[k] + P_{\text{batt}}[k] = 0 \label{hes_power_ERM}\\
P_{\text{batt}}[k] - P^{\text{c}}_{\text{batt}}[k] - P^{\text{d}}_{\text{batt}}[k] = 0\\
P^{\text{c}}_{\text{batt}}[k]~ P^{\text{d}}_{\text{batt}}[k] =0 \label{CS_ERM}\\
E[k] - E[k-1] - (1/E_{\text{C}})[\eta^{\text{batt}}_{\text{iv}}P^{\text{c}}_{\text{batt}}[k-1]\\
+ (1/\eta^{\text{batt}}_{\text{iv}})P^{\text{d}}_{\text{batt}}[k-1]] T_s=0\\
\underline{E} \le E[k] \le \bar {E} \label{SOC_ineq_ERM}\\
0\le P^{\text{c}}_{\text{batt}}[k] \le P^{\max}_{\text{batt}}\\
-P^{\max}_{\text{batt}}\le P^{\text{d}}_{\text{batt}}[k] \le 0\\
E[0]  = E_{\text{0}}, E[N] = E_{\text{0}} \label{sustain_ERM} 
\end{align}
\label{NLP_Temp_agn}
\end{subequations}
In this case, the temperature agnostic PV model, given by \eqref{P_cell_temp_agn} is used to calculate the $P_{\text{pv}}$ over the periodic horizon, $k\in\{1,\hdots,N\}$. The discretized PV equations are given as follows:
\begin{subequations}
\begin{align}
P^{\text{dc}}_{\text{pv}}[k]-({L~M~K}) \eta_{\text{cell}}~A_{\text{pv}}~S_{\text{I}}[k]=0\\
\eta^{\text{pv}}_{\text{iv}}[k]- \frac{\eta_{\text{max}}}{1+e^{-\gamma P^{\text{dc}}_{\text{pv}}[k]]}} = 0\\
P_{\text{pv}}[k]- \eta^{\text{pv}}_{\text{iv}}[k]P^{\text{dc}}_{\text{pv}}[k]=0 \label{PV2}
\end{align}
\label{PV_simplified}
\end{subequations}
\noindent Various case studies demonstrating the performance of different HES models with different battery ratings and under different ambient temperatures are presented and discussed in the next section.

\section{CASE STUDIES}\label{case_studies}
The codebase of the HES predictive optimization problem has been developed in Pyomo (a Python-based, open-source optimization modeling tool). The NLP problems (including projection, homotopy and temperature agnostic model) are solved using IPOPT (Interior Point OPTimizer) solver, while the simplified MIP problem has been solved using the Gurobi solver. The original NLP problem was warm-started by the solution of the simplified model. All the problems have been solved assuming constant ambient temperature throughout 24 hours. 

\subsection{With ISO New England energy price}
\begin{figure}[!htb]
\centering
\begin{subfigure}[!htb]{0.48\textwidth}
\centering
\includegraphics[width=\textwidth]{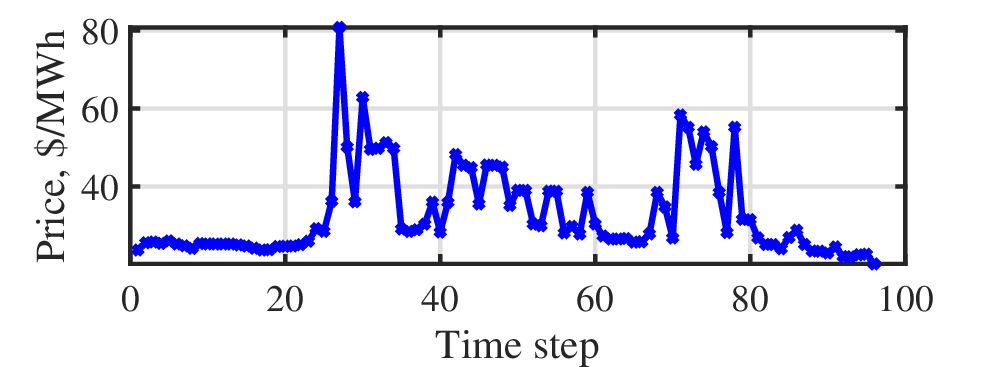}
\end{subfigure}

\begin{subfigure}[b]{0.48\textwidth}
\centering
\includegraphics[width=\textwidth]{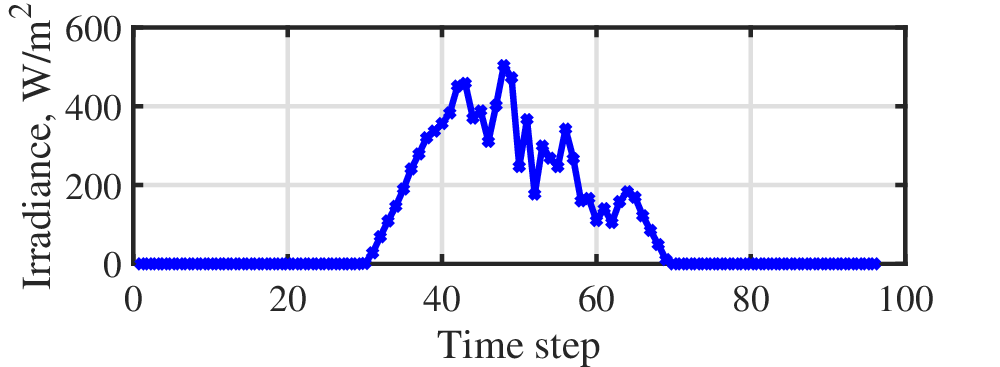}
\end{subfigure}
\caption{Input data: NE-ISO energy price and solar irradiance at 15-minute intervals for a day (24 hours).}
\label{Input data}
\end{figure}

The input data, $P^{r}$ for 15 minutes energy market is taken from ISO New England with a time step $T_{s}=0.25$ hour and $S_{\text{I}}$ are shown in Fig.\ref{Input data}. These data are only used to realize real-world settings and do not represent any particular event. Further, the HES's battery parameters are given in Table \ref{tab:battery parameters}. The initial values of state variable $SoC$, $T_{\text{batt}}$ and $T_{\text{en}}$ are assumed to be $E_{\text{0}}$, $T^{\text{0}}_{\text{batt}}$ and $T^{0}_{\text{en}}$ respectively. The PV consists of 120 modules and it is configured to form a 40 kW PV array. For temperature agnostic model, $\eta_{\text{cell}}$ of PV is $19.76\%$, $A_{\text{pv}}$=1.67 $m^2$ and energy capacity of the battery $E_{\text{c}}$ is taken as 43.2 kWh. To better align with real-world systems and observe the impact, we have scaled up our system by a factor of 1000, transitioning from a kW setup to a MW-scale configuration. 

\begin{table}[!htb]
\centering
\caption{Battery parameters}
\begin{tabular}{Sc Sc|Sc Sc} 
\toprule
 \thead{Parameter} & \thead{Value} & \thead{Parameter} & \thead{Value} \\
 \hline
 $i^{\text{max}}_{\text{batt}}$ & 50A (1C) & $E_{\text{0}}$ & 0.5\\

 $\bar{E}$ & 0.95 & $\underline{E}$ & 0.2\\

$\bar {v_{\text{batt}}}$ & 976 V &  $\underline {v_{\text{batt}}}$ & 714 V \\

 $T^{\text{min}}_{\text{batt}}$ & $15^{o}C$ &  $T^{\text{max}}_{\text{batt}}$& $35^{o}C$ \\

  $T^{\text{min}}_{\text{en}}$ & $15^{o}C$ &  $T^{\text{max}}_{\text{en}}$& $35^{o}C$ \\

 $C_{\text{Ah}}$ & 50 Ah & $C_{\text{T}}$  & 10 $KJ/^{o}C$ \\ 

$U$ & 0.2 $W/^{o}C$ & $U_{\text{en}}$& 0.001 $kW/^{o}C$ \\

 $R_{\text{0}}$ & 0.0716 $\Omega$ & $C_{\text{en}}$ & 30 $KJ/^{o}C$\\

 $N_{\text{cell}}$ & 100 & $T^{\text{0}}_{\text{batt}}$/ $T^{\text{0}}_{\text{en}}$& $20^{o}C$\\
\bottomrule
\end{tabular}
\label{tab:battery parameters}
\end{table}

\begin{figure}[!htb]
\centering
\includegraphics[width=0.48\textwidth]{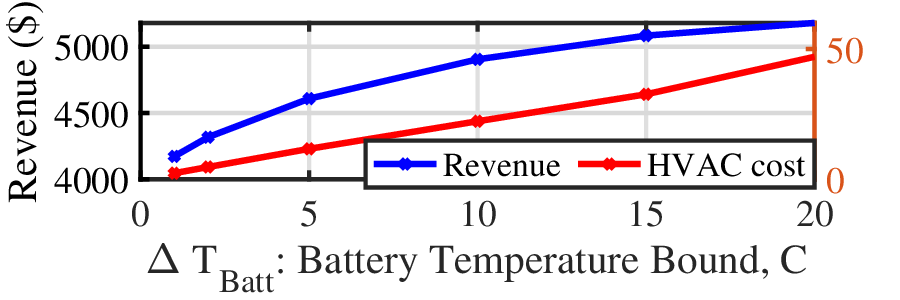}
\caption{Impact of variation in $\Delta T_{\text{Batt}}$ on overall revenue and HVAC cost at ambient temperature of $20^{o}$ C, 0.5C battery.}
\label{Temp bound variation}
\end{figure}

In the electro-thermal model, battery temperature plays a critical role while making the charging and discharging decisions as the $T_{\text{batt}}$ is correlated to $i_{\text{batt}}$ and $P_{\text{hv}}$ as shown in \eqref{Temp1} i.e. if the $T_{\text{batt}}$ hits the $T^{\text{max}}_{\text{batt}}$ or $T^{\text{min}}_{\text{batt}}$, the charging and discharging processes will be halted until the temperature falls back to a lower temperature level. Notably, the higher temperature limit will give more flexibility in charging and discharging but requires more HVAC power to maintain the temperature within the limits. The choice of an appropriate temperature limit is pivotal and hinges upon the objectives of predictive optimization. Fig.~\ref{Temp bound variation} shows the maximum HES revenue and expense of using HVAC by solving the NLP problem under varying $\Delta T_{\text{Batt}}$ condition. Here, $\Delta T_{\text{Batt}}$ represents the temperature bound i.e. $\Delta T_{\text{Batt}}$ = ($T^{\text{max}}_{\text{batt}}$~-~$T^{\text{min}}_{\text{batt}}$).

It is observed that the maximum revenue of the HES increases by 1.25 times while the HVAC cost increases by 20 times as the $\Delta T_{\text{Batt}}$ increases from $1^{0}C$ to $20^{0}C$. The increased limit will affect the battery's health by permitting rapid high charging and discharging, providing the option for complete discharges when prices are elevated. Thus, there exists a trade-off between revenue and battery health. In this paper,  values of $T^{\text{max}}_{\text{batt}}$ = $35^{0}C$ and $T^{\text{min}}_{\text{batt}}$= $15^{0}C$ have been selected.

\begin{figure}[!htb]
\centering
\begin{subfigure}[b]{0.5\textwidth}
\centering
\includegraphics[width=\textwidth]{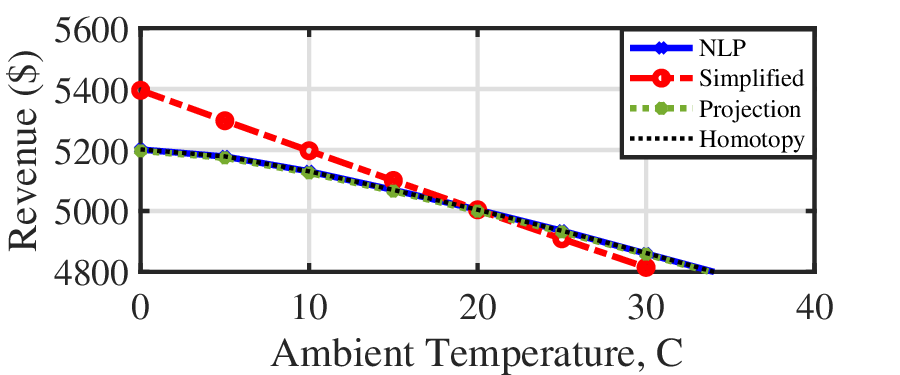}
\caption{0.25C battery }
\label{}
 \end{subfigure}
\begin{subfigure}[b]{0.5\textwidth}
\centering
\includegraphics[width=\textwidth]{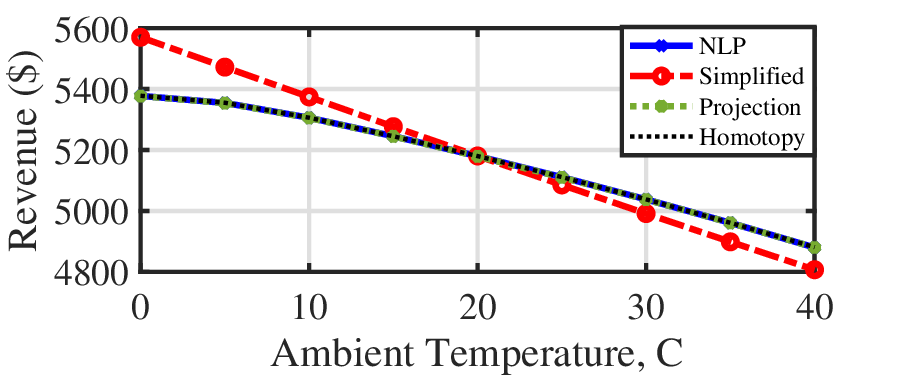}
\caption{0.5C battery}
\label{}
\end{subfigure}
\begin{subfigure}[b]{0.5\textwidth}
\centering
\includegraphics[width=\textwidth]{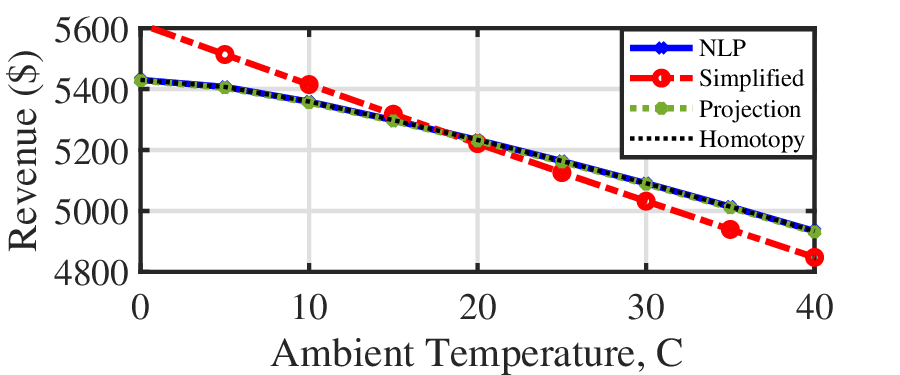}
\caption{1C battery }
\label{}
\end{subfigure}
\begin{subfigure}[b]{0.5\textwidth}
\centering
\includegraphics[width=\textwidth]{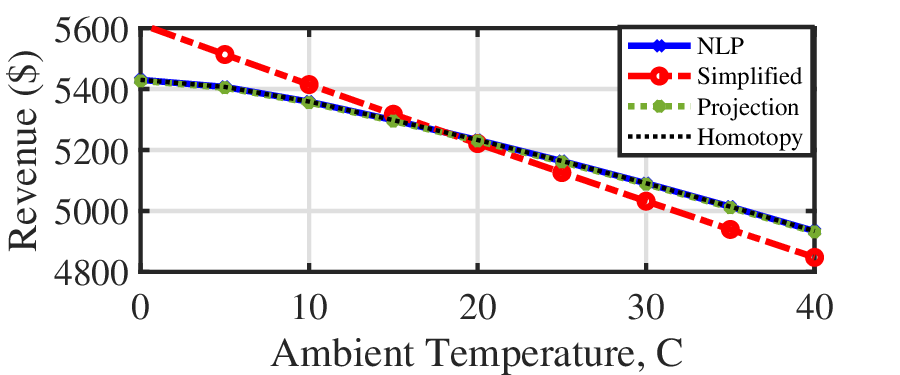}
\caption{4C battery}
\end{subfigure}
\caption{Predicted revenue comparison of various models under different ambient temperatures. Note that the simplified model's predicted revenue is not always realizable because it violates physical constraints.}
\label{Rev_comp1}
\end{figure}

The Fig.~\ref{Rev_comp1} shows the predicted revenue of the temperature-dependent optimization problem formulated in Section~\ref{Problem formulation}. The problem has been solved for different $T_{\text{amb}}$ and C-ratings of the battery.

It can be observed that the predicted revenue declines linearly for the simplified model and follows non-linear trends for NLPs as the $T_{\text{amb}}$ increases. 
The maximum deviation occurs at $T_{\text{amb}}$ = $0^{0}C$, when the simplified model predicts the $3.5\%$ extra revenue, which may not be the feasible schedule for HES and can lead to a penalty. This deviation is due to the simplified model's assumed constant $v_{\text{batt}}$. The solution obtained through projection has been close to the solution from the NLP solution. However, across all battery ratings, the NLP model consistently performs better. This shows that the NLP model provides a better solution within the feasible space. Furthermore, the homotopy demonstrates that it converges to the NLP solution, reinforcing the earlier claim that the NLP model finds the best solution.

\begin{figure}[!htb]
\centering
\begin{subfigure}[b]{0.48\textwidth}
\centering
\includegraphics[width=\textwidth]{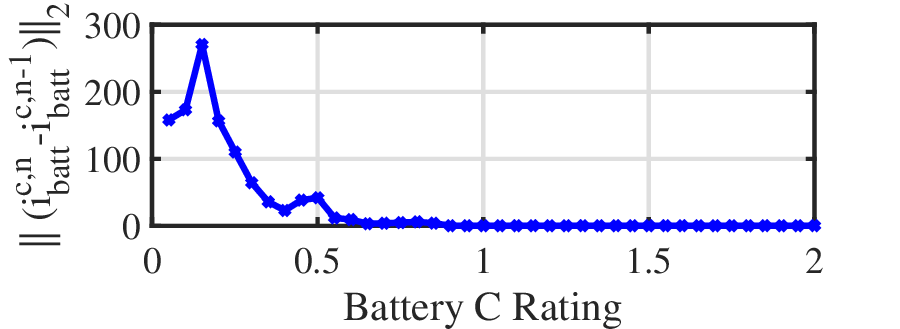}
\caption{decision variable: $i^{c}_{\text{batt}}$}
\label{}
 \end{subfigure}
 \hfill
\begin{subfigure}[b]{0.48\textwidth}
\centering
\includegraphics[width=\textwidth]{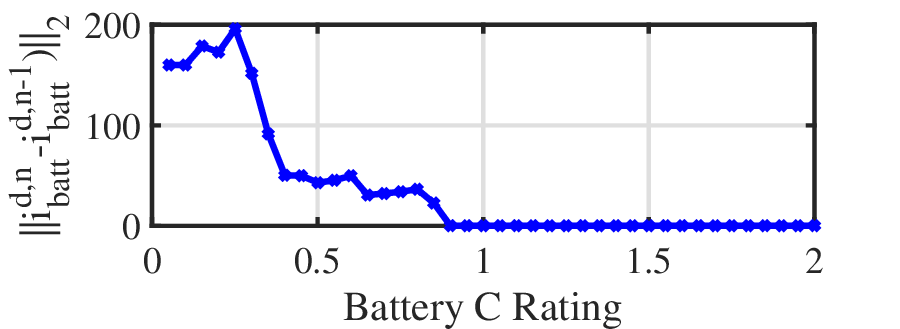}
\caption{decision variable: $i^{d}_{\text{batt}}$}
\label{}
\end{subfigure}
 \hfill
\begin{subfigure}[b]{0.48\textwidth}
\centering
\includegraphics[width=\textwidth]{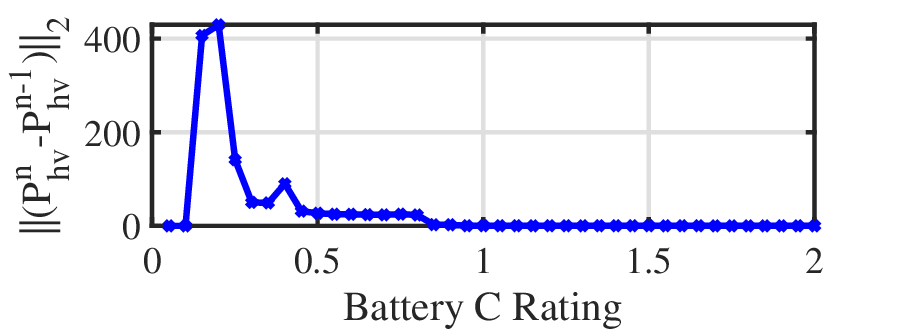}
\caption{decision variable: $P_{\text{hv}}$}
\label{}
\end{subfigure}
 \hfill
\begin{subfigure}[b]{0.48\textwidth}
\centering
\includegraphics[width=\textwidth]{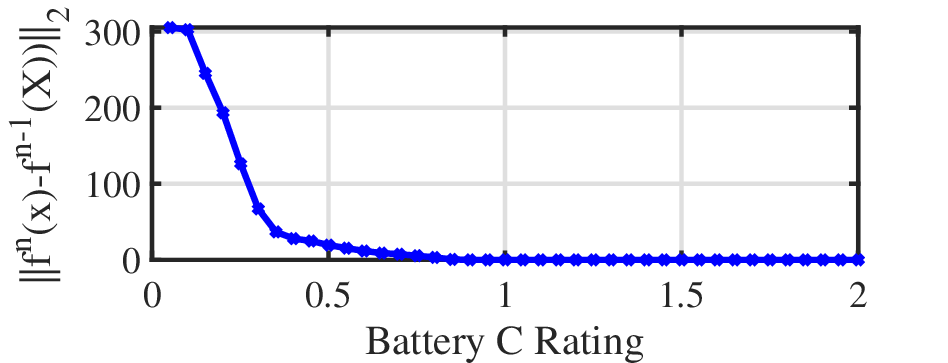}
\caption{Revenue}
\label{}
\end{subfigure}
\caption{Homotopy trajectory: 2-norm between previous and current solution on a normal day (at $T_{\text{amb}}$ = $20^{o}$ C).}
\label{homotopy_trajectory}
\end{figure}

The trajectory of the homotopy solution is shown in Fig.~\ref{homotopy_trajectory}. 

The 2-norm of the difference between the current (denoted by superscript 'n') and the previous (denoted by superscript 'n-1') decisions and the objective have been plotted. The figure is plotted until the  2C-rating for clarity as the 2-norm converges to 0, way before it. It can be observed that the solution is going from significant deviation to convergence to zero 2-norm, progressing towards the best solution by expanding the feasible solution space with each incremental increase in battery C-rating.

\begin{figure}[!htb]
\centering
\begin{subfigure}[b]{0.48\textwidth}
\centering
\includegraphics[width=\textwidth]{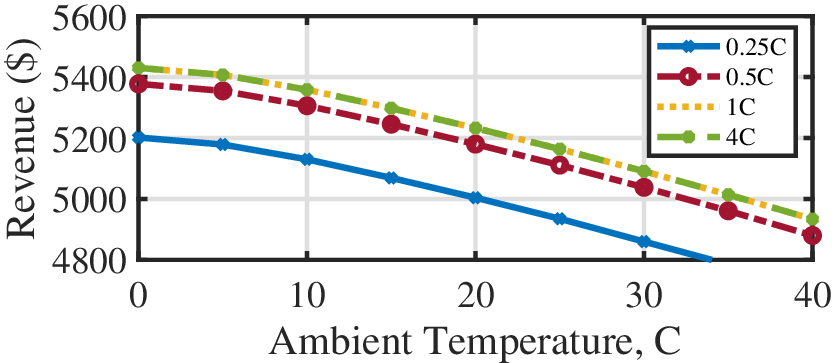}
\caption{NLP model}
\label{}
 \end{subfigure}
\begin{subfigure}[b]{0.48\textwidth}
\centering
\includegraphics[width=\textwidth]{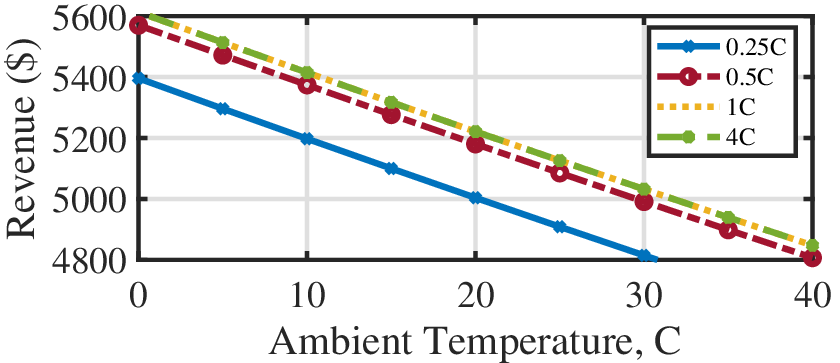}
\caption{Simplified model}
\label{}
\end{subfigure}
\begin{subfigure}[b]{0.48\textwidth}
\centering
\includegraphics[width=\textwidth]{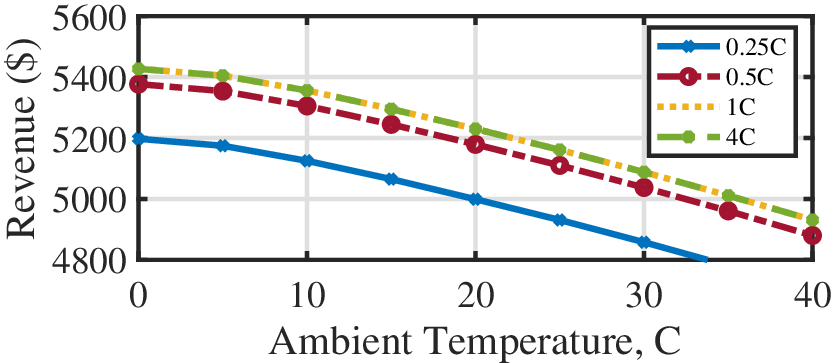}
\caption{Projection}
\label{}
\end{subfigure}
\begin{subfigure}[b]{0.48\textwidth}
\centering
\includegraphics[width=\textwidth]{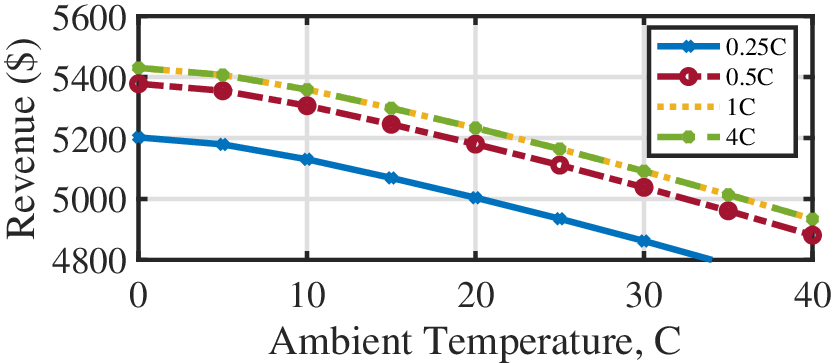}
\caption{Homotopy}
\end{subfigure}
\caption{Revenue comparison of different C-rated battery under different ambient temperatures. Note that the simplified model's predicted revenue is not always realizable because it violates physical constraints.}
\label{Rev_comp2}
\end{figure}

In Fig~\ref{Rev_comp2}, a progressive increase in predicted revenue has been observed with the transition from lower C-rated batteries to higher ones. 

Notably, the predicted revenue by HES employing 1C and 4C batteries closely align. This convergence is attributed to stringent temperature constraints; as these bounds become more restrictive, the higher charging speed associated with 4C batteries reaches a point of diminishing returns in revenue generation. Beyond a certain threshold, the faster charging does not significantly enhance the revenue due to thermal limitations.

The solutions from the temperature-agnostic models are very aggressive as they can charge or discharge very quickly ignoring the battery's thermal constraints. The predicted revenues of the temperature agnostic model for different C-rated batteries are given in Table~\ref{tab:Temp_agnostic_rev}. Here, The revenue increases with an increase in the C-rating of the battery which does not depict the charging limitations due to excessive heating of fast-charging batteries. The only feasible solution of the temperature agnostic model is with a 0.25C battery and if $T_{\text{amb}} \le 15^{o}C$, other schedules are not realizable considering the actual HES model.
\begin{table}[!htb]
\centering
\caption{Predicted revenue for temperature-agnostic model}
\begin{tabular}{Sc  Sc Sc|Sc Sc Sc} 
 \hline
 \thead{C-rating} & \thead{Revenue} & \thead{Feasible$^{\#}$} & \thead{C-rating} & \thead{Revenue}& \thead{Feasible$^{\#}$} \\
 \hline
 0.25C  & \$5054.62& $\checkmark$ & 1C & \$7261.83& $\times$\\
0.50C & \$5431.92& $\times$& 4C & \$8468.97 &$\times$\\
\hline
\multicolumn{6}{l}{$\#$  indicates whether this schedule is realizable with actual HES.}
\end{tabular}
\label{tab:Temp_agnostic_rev}
\end{table}

\subsection{Negative energy price scenario and PV curtailment}
\begin{figure}[!htb]
\centering
\includegraphics[width=0.4\textwidth]{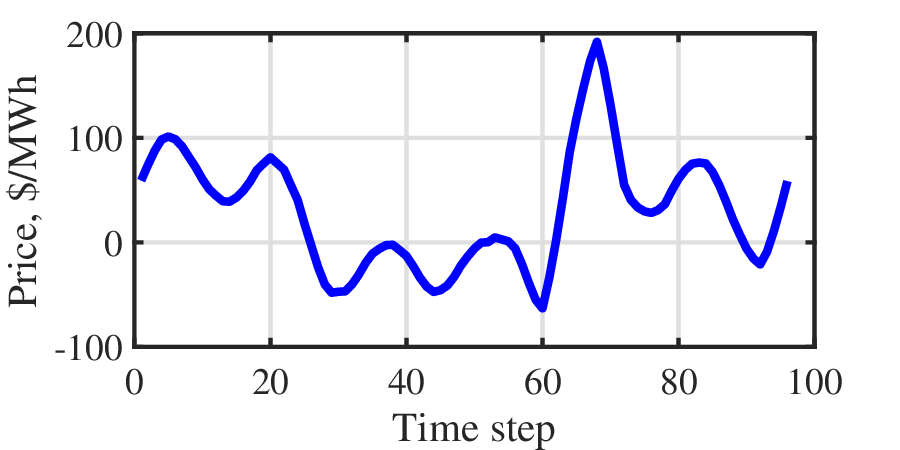}
\caption{Energy price }
\label{Case_2_Input}
\end{figure}

In the previous case, the energy prices were always positive. Hence, to test the robustness of the NLP and homotopy algorithm, test price $P^{r}$ has been used as shown in Fig.\ref{Case_2_Input}.

Negative energy prices can significantly reduce or even result in negative revenue for PV system owners. During such periods, PV generators might incur costs to offload excess electricity to grid operators or utilities. Curtailing PV generation in these circumstances is financially advantageous, preventing losses from selling surplus electricity at negative prices. Hence, another decision variable curtailment factor $\beta_{\text{cur}}$ is integrated in the optimization problem using \eqref{Curtailment}. The PV curtailment scheme is given by:
\begin{subequations}
\begin{align}
    P_{\text{pv}}[k]-\beta_{\text{cur}}[k] P^{\text{g}}_{\text{pv}}[k]  &=0\\
    0 \le \beta_{\text{cur}}[k] \le 1
\end{align}
\end{subequations}
where $P^{\text{g}}_{\text{pv}}$ denotes generated PV power, while $P_{\text{pv}}$ is total exported PV power to the grid. Table \ref{tab:Stress_case_rev_NO_cur} compares the revenue of all the temperature-dependent models and algorithms with and without PV curtailment under extreme temperatures and different C-ratings. It can be observed that the homotopy algorithm (highlighted using box) always converges to the same solution as the original NLP. Without the PV curtailment, the revenue increases as the ambient temperature increases due to lesser PV generation. 

\begin{table}[!htb]
\centering
\caption{Predicted revenue from the HES without PV curtailment } 
\begin{tabular}{Sc| Sc Sc| Sc Sc} 
\toprule

 \multirow{2}{*}{\thead{Method/Algorithm}} &  \multicolumn{2}{Sc|}{\thead{0.25C battery}} & \multicolumn{2}{Sc}{\thead{1C battery}} \\
\cline{2-5}
 & at $0^{o} C$ & at $40^{o} C$& at $0^{o} C$& at $40^{o} C$\\
\hline
NLP & \$5235.96 & \$5512.94 & \$5452.63 & \$5729.61\\
Homotopy &\fbox{\$5235.96} &\fbox{\$5512.94} & \fbox{\$5452.63}  & \fbox{\$5729.61}  \\
Simplified$^{*}$& \$5245.50 & \$5522.47 & \$5462.83  & \$5739.80\\
Projection &\$5231.01 & \$5507.99 &\$5448.43 & \$5725.41\\
%Temp Agnostic & \$165.52 & \$165.52& \$600.42 & \$600.42\\
\bottomrule
\multicolumn{5}{l}{*Simplified method's realized revenue is given by Projection value.}
\end{tabular}
\label{tab:Stress_case_rev_NO_cur}
\end{table}
Table \ref{tab:Stress_case_rev_cur} shows the predicted revenues with the PV curtailment. With the incorporation of the PV curtailment scheme, the predicted revenues grow under every condition as compared to those without curtailment. It can be observed from both the scenarios that simplified model always gives the upper bound of the solution but predicted revenues fall when the simplified solutions are projected in NLP space. 
\begin{table}[!htb]
\centering
\caption{Predicted revenue from the HES with PV curtailment } 
\begin{tabular}{Sc| Sc Sc| Sc Sc} 
\toprule
 \multirow{2}{*}{\thead{Method/Algorithm}} &  \multicolumn{2}{Sc|}{\thead{0.25C battery}} & \multicolumn{2}{Sc}{\thead{1C battery}} \\
\cline{2-5}
 & at $0^{o} C$ & at $40^{o} C$& at $0^{o} C$& at $40^{o} C$\\
\hline
NLP & \$7363.87 & \$7399.37 & \$7580.54 & \$7616.04\\
Homotopy &\fbox{\$7363.87} &\fbox{\$7399.37} & \fbox{\$7580.54}  & \fbox{\$7616.04}  \\
Simplified$^{*}$& \$7373.40 & \$7408.91 & \$7590.73  & \$7626.24\\
Projection &\$7358.89 & \$7394.39 &\$7576.31 & \$7611.81\\
%Temp Agnostic & \$165.52 & \$165.52& \$600.42 & \$600.42\\
\bottomrule
\multicolumn{5}{l}{*Simplified method's realized revenue is given by Projection value.}
\end{tabular}
\label{tab:Stress_case_rev_cur}
\end{table}

\subsection{Temperature agnostic vs Temperature dependent model}
\begin{figure}[!htb]
\centering
\includegraphics[width=0.45\textwidth]{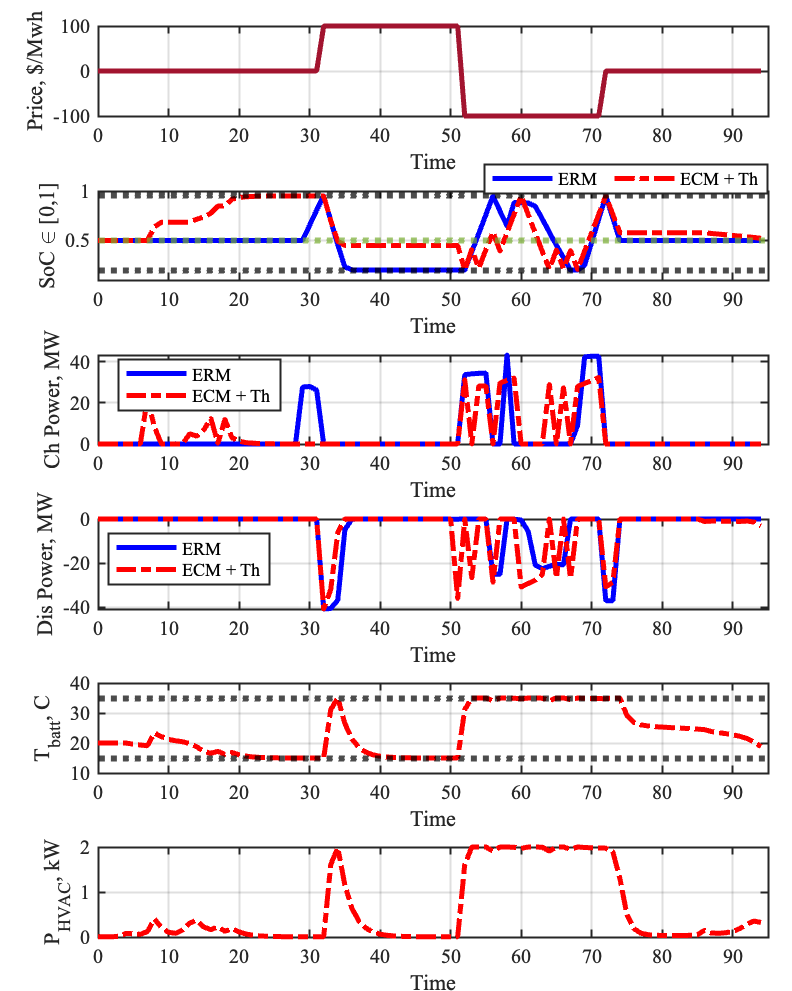}
\caption{Comparison of predicted scheduling of HES's 1C battery (ERM+Th) at $T_{\text{amb}}$ = 40$^{o}$C, U=5~W/$^{o}$C vs temperature-agnostic HES with ERM battery model. }
\label{U5_ERM}
\end{figure}

To assess the effectiveness of the temperature-dependent model for the HES, a comparison with the temperature-agnostic model (as discussed in Section \ref{temp_agnotic_section}) is crucial.
The temperature-agnostic model represents a temperature-agnostic PV+ battery model.
An energy price enforcing high discharging and charging rates with a 1C battery is employed. For this case study, objective function excludes considerations for HVAC usage, eliminating potential revenue loss. The thermal transmittance coefficient U represents diverse cooling systems, including air-cooled, surface-cooled, and liquid-cooled configurations. Lower coefficient values (U) signify slower battery cooling, while higher values represent faster cooling rates. The capital cost for HVAC systems for different cooling has been neglected.
This model is represented by ``ECM+Th". 
Table \ref{tab:TA vs TD} shows the predicted revenue for the temperature-agnostic and temperature-dependent model. For hot summer days ($T_{\text{amb}}$ = 40$^{o}$C), the temperature-agnostic model overestimates the revenue from PV by 12.5$\%$. Similar overestimation can be observed from the predicted battery revenue.
In Fig. \ref{U5_ERM}, a ECM+Th model with a high U value of 5 W/$^{o}$C is compared with the temperature-agnostic model (ERM). During fast cooling, generated heat dissipation within the enclosure necessitates higher HVAC power consumption to maintain the $T_{\text{en}}$ 
within specified limits. Consequently, with faster cooling rates, the battery can undergo rapid charging and discharging.
\begin{table}[t]
\centering
\caption{HES $\%$ change in predicted revenue comparing temperature-agnostic and temperature-dependent model with 1C battery, at $T_{\text{amb}}$ = 40$^{o}$C and U= 0.1 W/$^{o}$C, the capital costs of HVAC systems for different U has been neglected} 
\begin{tabular}{|Sc| Sc| Sc| Sc|} 
\toprule
\hline
 \thead{Model} &  \thead{HES (in  \%)}&\thead{PV (in \%)} & \thead{Battery (in \%)} \\
\hline
$\text{Temperature-agnostic}^{*}$& 14.5& 12.5& 16.8 \\
\hline
\hline
\multicolumn{4} {|Sc|}{Temperature-dependent}\\
\hline
U = 0.2 W/$^{o}$C& 0.2 & 0 & 0.42\\
\hline
U = 1 W/$^{o}$C & 1.25 & 0 & 2.67\\
\hline
U = 5 W/$^{o}$C & 2.6 & 0 & 5.45\\
\hline
\bottomrule
\multicolumn{4}{l}{*Not realizable with actual HES.}
\end{tabular}
\label{tab:TA vs TD}
\end{table}

% \begin{table}[!htb]
% \centering
% \caption{HES predicted revenue comparing temperature-agnostic and temperature-dependent model with 1C battery, at $T_{\text{amb}}$ = 40$^{o}$C, the capital costs of HVAC systems for different U has been neglected} 
% \begin{tabular}{|Sc| Sc| Sc| Sc|} 
% \toprule
% \hline
%  \thead{Model} &  \thead{HES (in  \$)}&\thead{PV (in \$)} & \thead{Battery (in \$)} \\
% \hline
% $\text{Temperature-agnostic}^{*}$& 14431& 7473& 6958 \\
% \hline
% \hline
% \multicolumn{4} {|Sc|}{Temperature-dependent}\\
% \hline
% U = 0.1 W/$^{o}$C&  12595 & 6640 & 5955\\
% \hline
% U = 0.2 W/$^{o}$C& 12620 & 6640 & 5980\\
% \hline
% U = 1 W/$^{o}$C & 12754 & 6640 & 6114\\
% \hline
% U = 5 W/$^{o}$C & 12920 & 6640 & 6280\\
% \hline
% \bottomrule
% \multicolumn{4}{l}{*Not realizable with actual HES.}
% \end{tabular}
% \label{tab:TA vs TD}
% \end{table}

\begin{figure}[!htb]
\centering
\includegraphics[width=0.48\textwidth]{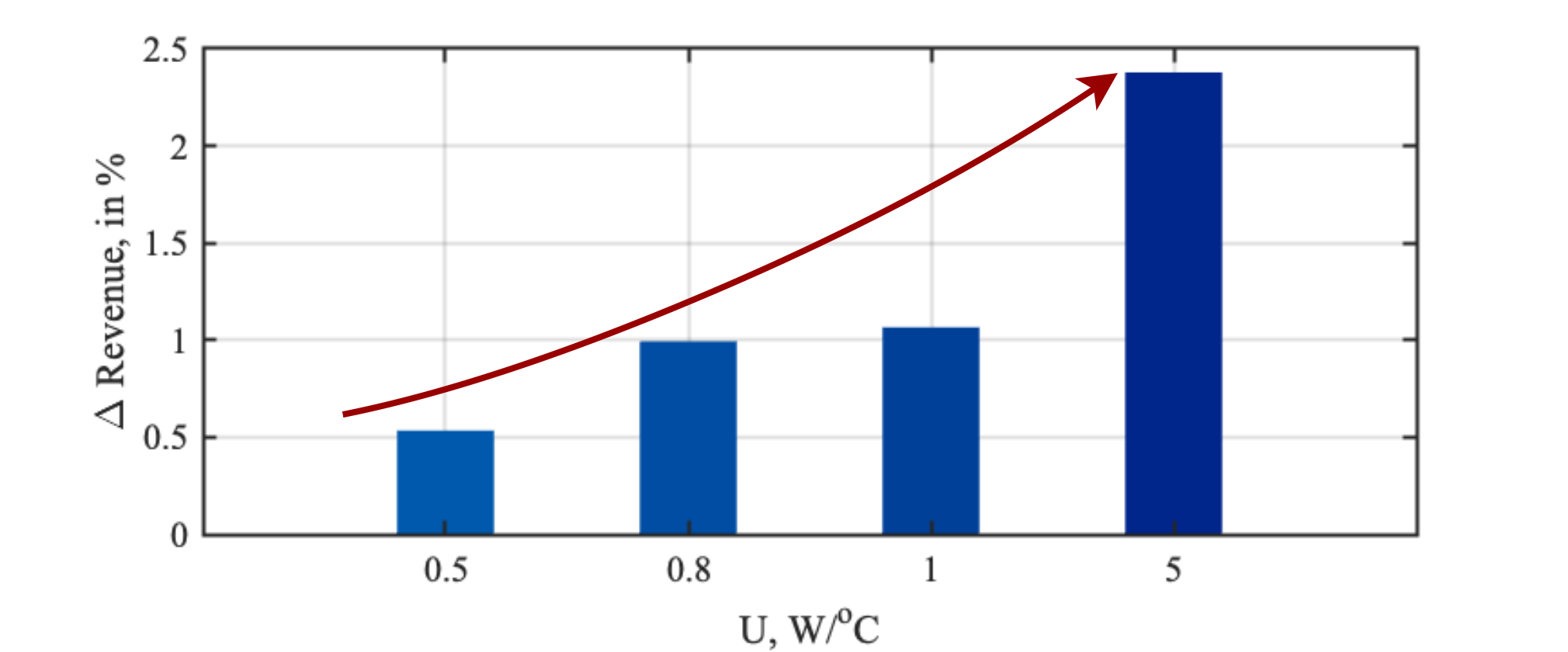}
\caption{Change in predicted revenue of temperature dependent model with HVAC is a free variable and variation of U, the base U= 0.2W/$^{o}$C. These results neglect any capital costs of HVAC systems for different U, which are expected to be non-trivial for $U\ge1$.}
\label{U_variation}
\end{figure}

As the $U$ value increases, the revenue predicted by the temperature-dependent model improves, as shown in Fig. \ref{U_variation}.
The change in revenue is calculated using the following equation: 
\begin{equation}
    \Delta Revenue~(\%) = \frac{Rev_{U} - Rev^{*}} {Rev^{*}} \times 100\\
\end{equation}
\noindent The $Rev^{*}$ represented the base value of predicted revenue when $U= 0.2$~W/$^{o}$C while $Rev_{U}$ represented the predicted revenue for $U \in$ [0.5, 0.8, 1, 5] W/$^{o}$C.
As the ERM model is temperature-agnostic, it operates without constraints imposed by thermal limits, resulting in overestimated predicted revenue that exceeds what can be realistically achieved by a physical battery. Fig. \ref{U5_ERM} shows that temperature-agnostic HES with ERM can discharge from $E_{\text{max}}$ to $E_{\text{min}}$ very quickly once the energy prices are suitable. At the same time, the temperature-dependent model is restricted by the $T_{\text{batt}}$. This led to an overestimation of predicted revenue by ERM, which the actual HES can't deliver. The performance of the temperature-dependent HES system is likely to be near the actual system. This indicates that utilities employing a temperature-agnostic model for predictive optimization may overestimate system performance, potentially leading to penalties for unmet expectations.

% \begin{figure}[!htb]
% \centering
% \includegraphics[width=1\columnwidth]{Overestimation.png}
% \caption{Temperature-agnostic model (ERM) exhibits performance overestimation compared to actual system ($U\in[0.1,5]$), while ECM+Th (U=0.5) close to the actual system.}
% \label{Overestimate}
% \end{figure}

\section{CONCLUSION}\label{Conclusion}
In this paper, a comprehensive electro-thermal model for HES has been thoroughly discussed and analyzed to give valuable insights. Within this framework, a predictive optimal scheduling problem to maximize HES revenue has been formulated. The intrinsic challenges arising from the system's nonlinear and non-convex nature, influenced by underlying physics and thermal constraints, make this problem inherently challenging. To address these complexities, a simplified problem using mixed-integer formulation has been proposed, which delivers the best solutions within defined constraints. The solution obtained from this simplified model is used to warm start the NLP. To validate the feasibility of the solution from the simplified model, a projection algorithm has been proposed. A homotopy algorithm has been proposed, which iterates over the C-rating of the battery. Lastly, to underscore the impact of temperature-dependent scheduling for HES, a temperature-agnostic model has been used. 

 These models are tested under various ambient temperatures, C-ratings and energy prices.  Based on the performance of the models, the homotopy solution displayed convergence towards zero 2-norm, as the feasible solution space gradually expands. The homotopy method also produces the same solution as the original NLP, while the simplified model provides the upper bound during stress testing. The predicted schedules from the temperature agnostic model are not realizable in actual plant settings.
 The discrepancy between the overestimated revenue projections of the temperature-agnostic model and the actual system performance underscores the importance of incorporating temperature-dependent. This work will be useful to give insight into HES complexities while highlighting the importance of a temperature-dependent model.

\bibliographystyle{IEEEtran}
\bibliography{PSCC_ref}

\end{document}